\documentstyle[epsf,12pt]{article}
\def\Bbb R{{\rm \bf R}}
\def\proclaim#1{\vskip2mm{\bf #1}\em}
\def\endproclaim{\em \vskip2mm}
\def\tag#1{\eqno(#1)}
\def\gathered{\begin{array}{c}}
\def\endgathered{\end{array}}
\def\text{\mbox}

\begin{document}

\title {A new formulation of the probe method\\ and related
problems}
\author{Masaru IKEHATA\\
Department of Mathematics,
Faculty of Engineering\\
Gunma University, Kiryu 376-8515, JAPAN}
\date{Final 24th December 2004}
\maketitle
\begin{abstract}
The probe method gives a general idea to obtain a reconstruction
formula of unknown objects embedded in a known background medium
from a mathematical counterpart (the Dirichlet-to-Neumann map) of
the measured data of some physical quantity on the boundary of the
medium. It is based on the sequence of special solutions of the
governing equation for the background medium related to a singular
solution of the equation.  In this paper the blowup property of
the sequence is clarified.  Moreover a new formulation of the probe
method based on the property is given in some typical inverse
boundary value problems.

\noindent
AMS: 35R30

\end{abstract}

\section{Introduction}

The probe method gives a general idea to obtain a reconstruction
formula of unknown objects embedded in a known background medium
from a mathematical counterpart (the Dirichlet-to-Neumann map) of
the measured data of some physical quantity on the boundary of the
medium. It was introduced by the author and applied to several
inverse boundary value problems and inverse scattering problems
(see \cite{I1,I2, I3, I4}).

The aim of this paper is to further investigate the probe method and
give a new formulation of the probe method,
which may be simpler than the previous formulation.
Since this paper is related to the idea of the probe method, we mainly consider only
a simple and typical inverse boundary value problem for the Helmholtz equation
which can be considered as a reduction of the inverse obstacle scattering problem, e.g.,
with point sources (see \cite{I3} for the reduction).

Let $\Omega$ be a bounded domain in $\Bbb R^m(m=2,3)$ with Lipschitz boundary.
Let $D$ be an open subset with Lipschitz boundary of $\Omega$ and
satisfy that $\overline D\subset\Omega$; $\Omega\setminus\overline D$ is connected.

We denote by $\nu$ the unit outward normal relative to $\Omega\setminus\overline D$.
Let $k\ge 0$.  We always assume that $0$ is not a Dirichlet eigenvalue of $\triangle +k^2$ in $\Omega$
and that $0$ is not an eigenvalue of the mixed problem
$$\begin{array}{c}
\displaystyle
\triangle u +k^2u=0\,\,\text{in}\,\Omega\setminus\overline D,\\
\\
\displaystyle
\frac{\partial u}{\partial\nu}=0\,\,\text{on}\,\,\partial D,\\
\\
\displaystyle
u=0\,\,\text{on}\,\partial\Omega.
\end{array}
$$
Given $f\in H^{1/2}(\partial\Omega)$ let $u\in H^1(\Omega\setminus\overline D)$ denote
the weak solution of the elliptic problem
$$\begin{array}{c}
\displaystyle
\triangle u+k^2u=0\,\,\text{in}\,\Omega\setminus\overline D,\\
\\
\displaystyle
\frac{\partial u}{\partial\nu}=0\,\,\text{on}\,\partial D,\\
\\
\displaystyle
u=f\,\,\text{on}\,\partial\Omega.
\end{array}
$$
Define
$$\displaystyle
\Lambda_Df=\frac{\partial u}{\partial\nu}\vert_{\partial\Omega}.
$$
We set $\Lambda_D=\Lambda_0$ in the case when $D=\emptyset$.
\noindent
$\Lambda_D$ is called the Dirichlet-to-Neumann map.

Here we consider the problem of extracting information about the shape and location of $D$
from $\Lambda_D$ or its partial knowledge.
The probe method gives us a reconstruction formula of
$\partial D$ by using $(\Lambda_0-\Lambda_D)f$ for infinitely many $f$.

For the description we need two concepts: needle and impact parameter.
A continuous curve $c:[0,\,1]\longmapsto\overline\Omega$ is called a needle
if $c(0), c(1)\in\partial\Omega$ and $c(t)\in\Omega$ for all $t\in\,]0,\,1[$.
Set $c_t=\{c(s)\,\vert\,0<s\le t\}$.
Define the impact parameter of $c$ with respect to $D$ by the formula
$$
t(c;D)=\sup\{0<t<1\,\vert\,\forall s\in\,]0,\,t[\,c(s)\in\Omega\setminus\overline D\}.
$$
If $t(c;D)<1$, then the impact parameter coincides with the first hitting parameter
of the curve $c$ with respect to $D$; if $t(c;D)=1$, then this means that the curve $c([0,\,1])$
is outside $D$.

We denote by $G_k(x)$ the standard fundamental solution of the Helmholtz equation.

\noindent
The starting point is the following.

\proclaim{\noindent Proposition 1.1.}
Given a needle $c$ and $t\in\,]0,\,1[$ there exists a sequence $v_1(\,\cdot\,;c_t), v_2(\,\cdot\,;c_t),\cdots$
of $H^1(\Omega)$ solutions of the Helmholtz equation such that, for each fixed
compact set $K$ of $\Bbb R^m$ with $K\subset\Omega\setminus c_t$
$$\displaystyle
\lim_{n\longrightarrow\infty}(\Vert v_n(\,\cdot\,;c_t)-G_k(\,\cdot\,-c(t))\Vert_{L^2(K)}
+\Vert\nabla\{v_n(\,\cdot\,;c_t)-G_k(\,\cdot\,-c(t))\}\Vert_{L^2(K)})=0.
$$
\endproclaim

\noindent
This is a consequence of Theorem 4 in \cite{I3}
which states the Runge approximation property
for the stationary Schr\"odinger equation (see also appendix A.1).

Define
$$\displaystyle
I_n(t;c)=\int_{\partial\Omega}\{(\Lambda_0-\Lambda_D)\overline f_n\}f_ndS
\tag {1.1}
$$
where
$$
f_n(y)=v_n(y;c_t),\,\,y\in\partial\Omega.
$$
We write
$$\displaystyle
I(t;c)=\lim_{n\longrightarrow\infty}I_n(t;c)
$$
if it exists.  This is called the indicator function.

\noindent
Define
$$\displaystyle
T(c)=\{t\in\,]0,\,1[\,\vert\, \forall\,s\in\,]0,\,t[\,
I(s;c)\,\text{exists and}\,\sup_{0<s<t}I(s;c)<\infty\}.
\tag {1.2}
$$
We have already established the following \cite{I3}.

\proclaim{\noindent Theorem 1.1.}
Assume that both $\partial\Omega$ and $\partial D$ are $C^2$.
Then, for any needle $c$ the formula
$$\displaystyle
T(c)=]0\,,t(c;D)[,
\tag {1.3}
$$
is valid.

\endproclaim

\noindent
Since we have the formula
$$
\partial D=\{c(t)\vert_{t=t(c;D)}\,\vert\,t(c;D)<1\},
$$
we obtain the reconstruction formula of $\partial D$
from $\Lambda_D$ through (1.1), (1.2) and (1.3).

\noindent
This is the original formulation of the result obtained by applying
the probe method.  From this theorem we know that
$I(t;c)=\lim_{n\longrightarrow\infty}I_n(t;c)$ exists if
$0<t<t(c;D)$. In addition, it is easy to see that $\lim_{t\uparrow
t(c;D)}I(t;c)=\infty$ in the case when $t(c;D)<1$. However, if $1>t\ge t(c;D)$, we did not
mention explicitly the behaviour of $I_n(t;c)$ as
$n\longrightarrow\infty$ in the papers devoted to the probe
method.

Recently Erhard-Potthast \cite{EP} studied the probe method
numerically. They considered, as an example, an inverse boundary
value problem for the Helmholtz equation for sound-soft obstacles
($u=0$ on $\partial D$) and computed an approximation of the
corresponding indicator function by employing the techniques of the
point source and singular sources methods by Potthast \cite{P1,P2}.
Their computation results show that the absolute value of the
approximation takes a large value when $t>t(c;D)$ and $c(t)\in D$.
This suggests the blowup of the indicator function when the
parameter $t$ in the indicator function is greater than the impact
parameter and the corresponding point on the needle inside the
unknown objects.

\begin{figure}[htbp]
\begin{center}
\epsfxsize=10cm
\epsfysize=14cm
\epsfbox{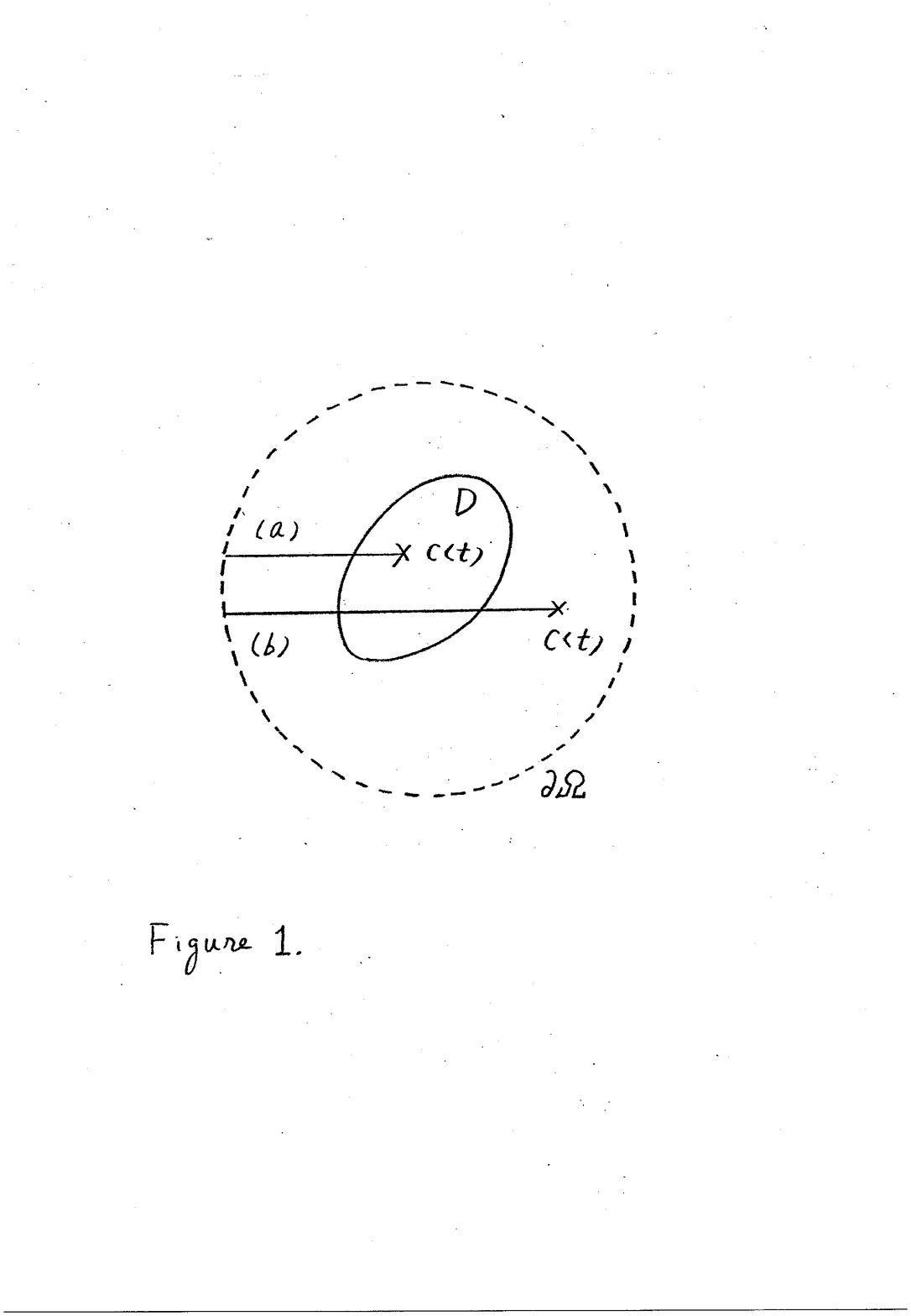}
\caption{(a) $c(t)\in D$.  (b) $c(t)\in\Omega\setminus\overline D$.
Both cases satisfy $t(c;D)<t<1$.
}\label{fig1}
\end{center}
\end{figure}

In this paper we give the proof of the blowup property of the
indicator function provided $k$ is small enough. More precisely,
we obtain: if $t(c;D)<1$ and $1>t\ge t(c;D)$, then
$\lim_{n\longrightarrow\infty}I_n(t;c)=\infty$ under suitable
conditions on $c$. If $c(t)\in D$, then this result gives a
verification of Erhard-Potthast's computation result.  However,
our result covers the case also when $c(t)$ is outside $D$
(see Figure \ref{fig1} for the geometry).

\newpage


\vspace{-1.0cm}

\noindent
This is an unexpected property of the indicator function and needs
purely theoretical consideration. Their
computation result does not cover this case since their approximation of the indicator
function is too simple. The result is based on the discovery
of the blowup of the sequence of the solutions of the Helmholtz
equation given in Proposition 1.1 on the needle (Lemmas 2.1 and
2.2).

However, for the description of the result we do not make use of the formulation given above.
We give a new and simpler formulation of the probe method.
In the formulation, we do not make use of the impact parameter.

\section{New formulation of the probe method}

In this section, we introduce a new formulation of the probe
method.  Given a point $x\in\Omega$ let $N_x$ denote the set of all piecewise linear
curves $\sigma:[0,\,1]\longmapsto\overline\Omega$ such
that

(1) $\sigma(0)\in\partial\Omega$, $\sigma(1)=x$ and
$\sigma(t)\in\Omega$ for all $t\in]0,\,1[$;

(2) $\sigma$ is injective.

\noindent
We call $\sigma\in N_x$ a needle with tip at $x$.
For the new formulation of the probe method we need the following.

{\bf\noindent Definition 2.1.}
Let $\sigma\in N_x$.  We call the sequence $\xi=\{v_n\}$
of $H^1(\Omega)$ solutions of the Helmholtz equation
a needle sequence for $(x,\sigma)$ if it satisfies,
for each
fixed compact set $K$ of $\Bbb R^m$ with $K\subset\Omega\setminus
\sigma([0,\,1])$
$$\displaystyle
\lim_{n\longrightarrow\infty}(\Vert v_n(\,\cdot\,)-G_k(\,\cdot\,-x)\Vert_{L^2(K)}
+\Vert\nabla\{v_n(\,\cdot\,)-G_k(\,\cdot\,-x)\}\Vert_{L^2(K)})=0.
$$
\noindent
Needless to say, the existence of the needle sequence is a consequence of
Proposition 1.1.  The problem is the behaviour of the needle sequence
on the needle as $n\longrightarrow\infty$.

Here we make a definition.
Let $\mbox{\boldmath $b$}$ be a nonzero vector in $\Bbb R^m$.
Given $x\in\Bbb R^m$, $\rho>0$ and $\theta\in]0,\pi[$
the set
$$\displaystyle
V=\{y\in\Bbb R^m\,\vert\,\vert y-x\vert<\rho\,\,\text{and}\,\,
(y-x)\cdot\mbox{\boldmath $b$}>\vert y-x\vert\vert\mbox{\boldmath $b$}\vert
\cos(\theta/2)\}
$$
is called a finite cone of height $\rho$, axis direction $\mbox{\boldmath $b$}$ and
aperture angle $\theta$ with vertex at $x$.

\noindent
The two lemmas given below are the core of the new formulation
of the probe method.

\proclaim{\noindent Lemma 2.1.}
Let $x\in\Omega$ be an arbitrary point and $\sigma$ be a needle with tip at $x$.
Let $\xi=\{v_n\}$ be an arbitrary needle sequence for $(x,\sigma)$.
Then, for any finite cone $V$ with vertex at $x$ we have
$$\displaystyle
\lim_{n\longrightarrow\infty}\int_{V\cap\Omega}\vert\nabla v_n(y)\vert^2dy=\infty.
$$
\endproclaim

{\it\noindent Proof.}
We employ a contradiction argument.
Assume that the conclusion is not true.
Then there exist $M>0$ and a sequence $n_1<n_2<\cdots\longrightarrow\infty$ such that
$$\displaystyle
\int_{V\cap\Omega}\vert\nabla v_{n_j}(y)\vert^2dy<M, j=1,2,\cdots.
$$
Take a sufficiently small open ball $B$ centred at $x$ with radius $R$ such that
$\overline B\subset\Omega$ and $\sigma(]0,\,1])\cap B$ becomes a segment having $x$ as an
end point.  Then one can find a finite cone $V'\subset V$ with vertex at $x$ such that,
for every $\epsilon$ with $0<\epsilon<R$ $K_{\epsilon}\equiv\overline V'\cap(\overline B\setminus B_{x}(\epsilon))
\subset V\cap(\Omega\setminus\sigma(]0,\,1]))$ where $B_x(\epsilon)$ stands for the open ball centred
at $x$ with radius $\epsilon$.
Thus we have
$$\displaystyle
\int_{K_{\epsilon}}\vert\nabla v_{n_j}(y)\vert^2dy<M, j=1,2,\cdots.
$$
Since $\nabla v_{n_j}(\,\cdot\,)\longrightarrow\nabla G_k(\,\cdot\,-x)$ in $L^2(K_{\epsilon})$,
we get
$$\displaystyle
\int_{K_{\epsilon}}\vert\nabla G_k(y-x)\vert^2 dy\le M.
$$
Since $\epsilon$ can be arbitrary small, applying Fatou's lemma for $\epsilon=1/l$ as $l\longrightarrow\infty$
to the integral, we obtain
$$\displaystyle
\int_{V'\cap B}\vert\nabla G_k(y-x)\vert^2 dy\le M.
$$
However, using polar coordinates centred at $x$ one can show that this left hand side is divergent.
This is a contradiction and completes the proof.

\noindent
$\Box$

\proclaim{\noindent Lemma 2.2.}
Let $x\in\Omega$ be an arbitrary point and $\sigma$ be a needle with tip at $x$.
Let $\xi=\{v_n\}$ be an arbitrary needle sequence for $(x,\sigma)$.
Then for any point $z\in\sigma(]0,\,1[)$ and open ball $B$ centred at $z$ we have
$$\displaystyle
\lim_{n\longrightarrow\infty}\int_{B\cap\Omega}\vert\nabla v_n(y)\vert^2dy=\infty.
$$
\endproclaim

{\it\noindent Proof.}

\noindent Let $v$ be an
arbitrary solution of the Helmholtz equation in $\Omega$.  Note
that $v$ can be identified with a smooth function in $\Omega$ and
all the derivatives satisfy the Helmholtz equation in $\Omega$.
Choose an open ball $B'$ centred at $z$ such that $\overline
B'\subset B\cap\Omega$. Next choose a smaller open ball $B''$ centred at $z$ such that
$\overline B''\subset B'$.  Applying (A.1) to the case when $W=B'$
and $K=\overline B''$, we have
$$\displaystyle
\int_{B''}\vert\nabla(\nabla v)\vert^2dy\le C\int_{B'}\vert\nabla v\vert^2 dy.
\tag {2.1}
$$
Applying the trace theorem to $B''$, we have
$$\displaystyle
\int_{\partial B''}\vert\nabla v\vert^2 dS
\le C'(\int_{B''}\vert\nabla v\vert^2 dy
+\int_{B''}\vert\nabla(\nabla v)\vert^2 dy).
\tag {2.2}
$$
Choose a $C^2$ domain $U$ in such a way that $\Sigma\equiv\partial
U\cap\partial B''$ has a positive surface measure on $\partial B''$,
$\text{dist}\,(\partial
U\setminus\Sigma, \sigma)>0$, $x\in U$ and $\vert U\vert$ is sufficiently small
in the following sense:
$$\displaystyle
\omega_m>k^m\vert U\vert
$$
where $\omega_m$ is the volume of the unit ball in $\Bbb R^m$.
This last inequality implies that $0$ is not a Dirichlet eigenvalue
of $\triangle+k^2$ in $U$ (see Lemma 1 in \cite{SU}).

\begin{figure}[htbp]
\begin{center}
\epsfxsize=10cm
\epsfysize=14cm
\epsfbox{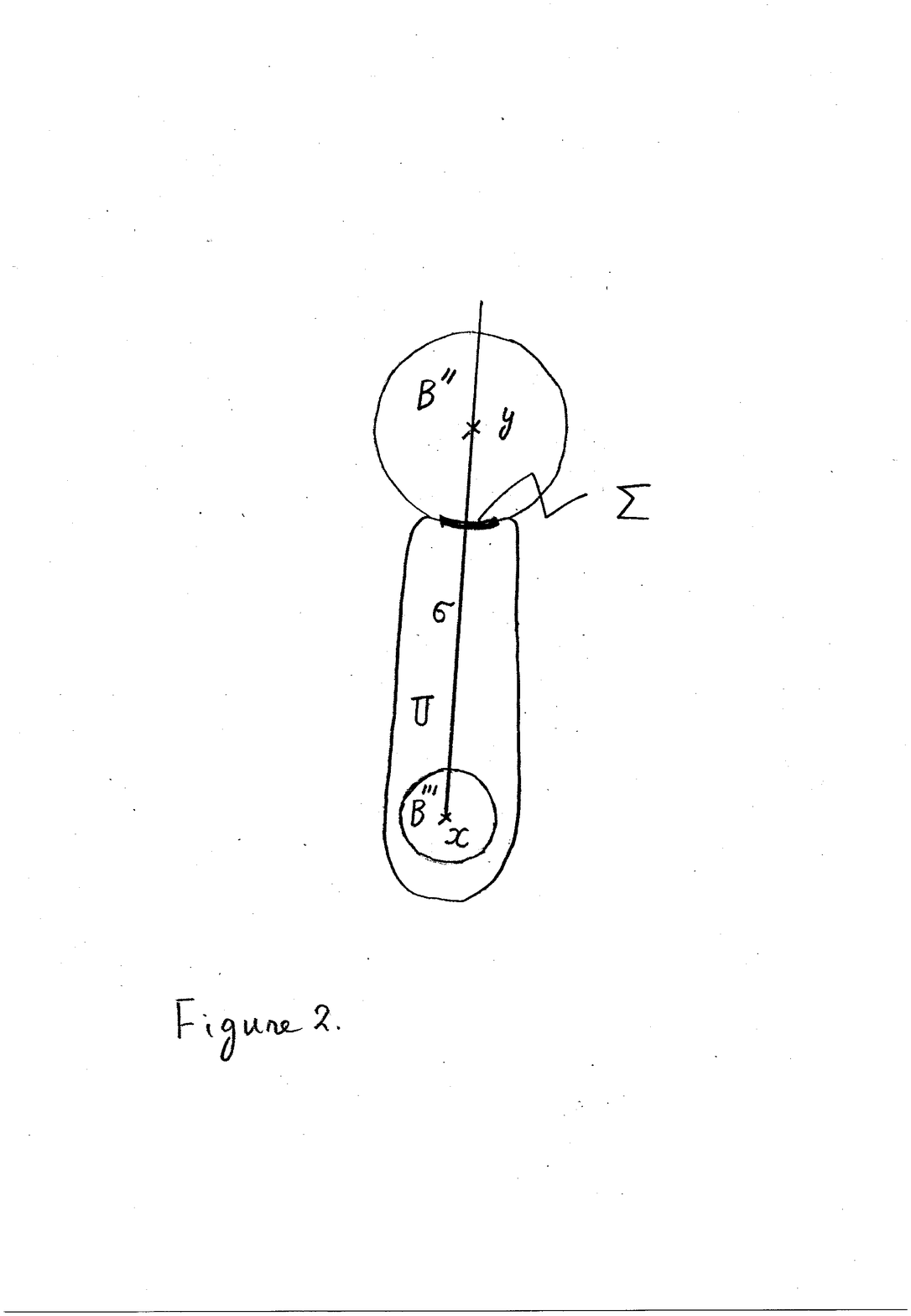}
\caption{
An illustration of $B'', U, \Sigma$ and $B'''$.
}\label{fig2}
\end{center}
\end{figure}

Choose an open ball $B'''$ centered
at $x$ such that $\overline B'''\subset U$
(see Figure \ref{fig2} for the geometry).

\newpage

\vspace{-1.0cm}

\noindent
Then (A.2) for the
case when $W=U$ gives
$$\begin{array}{c}
\displaystyle
\int_{B'''}\vert\nabla v\vert^2dy\\
\\
\displaystyle
\le\int_U\vert\nabla v\vert^2dy\\
\\
\displaystyle
\le C''\int_{\partial U}\vert\nabla v\vert^2dS\\
\\
\displaystyle
=C''(\int_{\Sigma}\vert\nabla v\vert^2dS+\int_{\partial U\setminus\Sigma}\vert\nabla v\vert^2dS)
\\
\\
\displaystyle
\le C''(\int_{\partial B''}\vert\nabla v\vert^2dS
+\int_{\partial U\setminus\Sigma}\vert\nabla v\vert^2dS).
\end{array}
\tag {2.3}
$$
From (2.1), (2.2) and (2.3) we obtain the estimate of $\nabla v$ in $B'$ in terms of $\nabla v$
in $B'''$ from below:
$$\displaystyle
\int_{B'''}\vert\nabla v\vert^2dy
\le C'''(\int_{\partial U\setminus\Sigma}\vert\nabla v\vert^2dS
+\int_{B'}\vert\nabla v\vert^2 dy).
\tag {2.4}
$$
Now set $v=v_n(\,\cdot\,)$.  Since $\text{dis}\,(\partial U\setminus\Sigma,\sigma(]0,\,1]))>0$
and $\nabla v_n(\,\cdot\,)$ converges to $\nabla G_k(\,\cdot\,-x)$
in $H^1_{\text{loc}}(\Omega\setminus\sigma(]0,\,1]))$,
the trace theorem gives
$$
\lim_{n\longrightarrow\infty}\int_{\partial U\setminus\Sigma}\vert\nabla v_n(y)\vert^2dy
=\int_{\partial U\setminus\Sigma}\vert\nabla G_k(y-x)\vert^2 dy<\infty.
\tag {2.5}
$$
On the other hand, from Lemma 2.1, one knows that
$$\displaystyle
\lim_{n\longrightarrow\infty}\int_{B'''}\vert\nabla v_n(y)\vert^2dy=\infty.
\tag {2.6}
$$
Thus from (2.4) for $v=v_n(\,\cdot\,)$, (2.5) and (2.6) we obtain the desired conclusion.

\noindent
$\Box$

The argument given above can be applied to other elliptic equations and the elliptic systems by a
suitable modification.

A combination of Lemmas 2.1 and 2.2 tells us that any
needle sequence for any needle blows up on the needle.  The needle
sequence behaves like a beam!  This is a new fact not mentioned in
the previous papers about the probe method.

In order to describe our main result we introduce two positive constants
appearing in two types of the Poincar\'e inequalities (e.g., see \cite{Z}).
One is given in the following.

\proclaim{\noindent Proposition 2.1.}
For all $w\in H^1(\Omega\setminus\overline D)$ with $w=0$ on $\partial\Omega$
$$
\int_{\Omega\setminus\overline D}\vert w\vert^2 dy\le C_0(\Omega\setminus\overline D)^2
\int_{\Omega\setminus\overline D}\vert\nabla w\vert^2 dy
$$
where $C_0(\Omega\setminus\overline D)$ is a positive constant independent of $w$.
\endproclaim

{\it\noindent Proof.}  This is nothing but a standard compactness argument.

\noindent
$\Box$

The dependence of $C_0(\Omega\setminus\overline D)$ on $\Omega\setminus\overline D$
should be clarified.  However it is not the aim of this paper.
Another is given in the following.

\proclaim{\noindent Proposition 2.2.} Let $U$ be a bounded Lipschitz domain
of $\Bbb R^m$. For any $v\in H^1(U)$ we have
$$
\displaystyle
\int_U\vert v-v_U\vert^2 dy\le C(U)^2\int_U\vert\nabla v\vert^2 dy
$$
where $C(U)$ is a positive constant independent of $v$
and
$$
\displaystyle
v_U=\frac{1}{\vert U\vert}\int_U vdy.
$$

\endproclaim
{\it\noindent Proof.}
Again, this is nothing but a standard compactness argument.

\noindent
$\Box$

As a corollary we have

\proclaim{\noindent Proposition 2.3.} Let $U$ be a bounded Lipschitz domain
of $\Bbb R^m$. For any $v\in H^1(U)$ and Lebesgue
measurable $A\subset U$ with $\vert A\vert>0$ we have
$$
\displaystyle
\int_U\vert v-v_A\vert^2 dy\le C(U)^2(1+\frac{\vert U\vert^{1/2}}{\vert A\vert^{1/2}})^2\int_U\vert\nabla v\vert^2 dy
$$
where $C(U)$ is the same constant as that of Proposition 2.2
and
$$
\displaystyle
v_A=\frac{1}{\vert A\vert}\int_A vdy.
$$
\endproclaim
{\it\noindent Proof.}  The following argument is taken from \cite{SS}(see also \cite{Z} for an abstract version).
Proposition 2.2 gives
$$\begin{array}{c}
\displaystyle
\Vert v-v_A\Vert_{L^2(U)}
\le\Vert v-v_U\Vert_{L^2(U)}
+\Vert\frac{1}{\vert A\vert}\int_{A}(v-v_U)dy\Vert_{L^2(U)}\\
\\
\displaystyle
\le
C(U)\Vert\nabla v\Vert_{L^2(U)}
+\frac{\vert U\vert^{1/2}}{\vert A\vert}
\vert\int_A(v-v_U)dy\vert\\
\\
\displaystyle
\le
C(U)\Vert\nabla v\Vert_{L^2(U)}
+\frac{\vert U\vert^{1/2}}{\vert A\vert^{1/2}}
\Vert v-v_U\Vert_{L^2(U)}.
\end{array}
$$
Then again Proposition 2.2 gives the desired estimate.

\noindent
$\Box$

We make use of the property that $\displaystyle C(U)^2(1+\frac{\vert U\vert^{1/2}}{\vert A\vert^{1/2}})^2$
continuously depends on $\vert A\vert$ for each fixed $U$.

\noindent
{\bf\noindent Definition 2.2.}
Given $x\in\,\Omega$, needle $\sigma$ with tip $x$
and needle sequence $\xi=\{v_n\}$ for $(x,\sigma)$
define
$$\displaystyle
I(x,\sigma,\xi)_n=\int_{\partial\Omega}\{(\Lambda_0-\Lambda_D)\overline f_n\}
f_ndS,\,n=1,2,\cdots
$$
where
$$\displaystyle
f_n(y)=v_n(y),\,\,y\in\,\partial\Omega.
$$

\noindent $\{I(x,\sigma,\xi)_n\}_{n=1,2,\cdots}$ is a sequence
depending on $\xi$ and $\sigma\in N_x$.
We call the sequence the indicator sequence.

Now the main result is the following.

\proclaim{\noindent Theorem 2.1.}
Assume that $D$ is given by a union of finitely many bounded Lipschitz domains
$D_1,\cdots, D_N$ such that $\overline D_j\cap\overline D_l=\emptyset$ if $j\not=l$.
Let $k\ge 0$ be small in the following sense:
$$
k^2C_0(\Omega\setminus\overline D)^2\le 1
\tag {2.7}
$$
and
$$\min_j\{1-2k^2C(D_j)^2(1+1)^2\}>0.
\tag {2.8}
$$
Then, given $x\in\Omega$ and needle $\sigma$ with tip at $x$ we have:

if $x\in\Omega\setminus\overline D$ and $\sigma(]0,1])\cap\overline D=\emptyset$,
then for any needle sequence $\xi=\{v_n\}$ for $(x,\sigma)$
the sequence $\{I(x,\sigma,\xi)_n\}$ is convergent;

if $x\in\Omega\setminus\overline D$ and $\sigma(]0,1])\cap D\not=\emptyset$,
then for any needle sequence $\xi=\{v_n\}$ for $(x,\sigma)$ we have
$\lim_{n\longrightarrow\infty}I(x,\sigma,\xi)_n=\infty$;

if $x\in\overline D$,
then for any needle sequence $\xi=\{v_n\}$ for $(x,\sigma)$ we have
$\lim_{n\longrightarrow\infty}I(x,\sigma,\xi)_n=\infty$.

\endproclaim

\begin{figure}[htbp]
\begin{center}
\epsfxsize=10cm
\epsfysize=14cm
\epsfbox{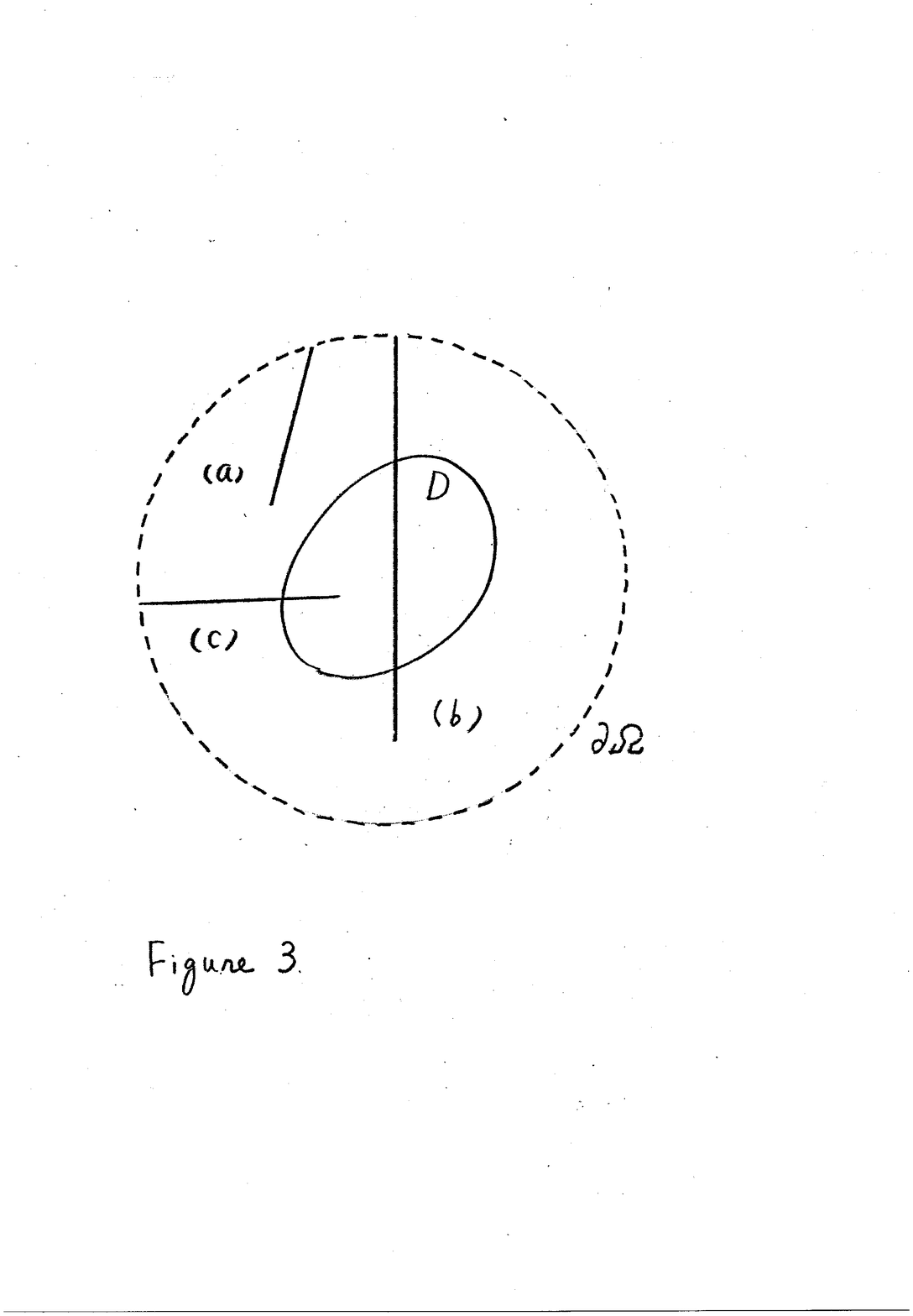}
\caption{
 An illustration of three cases: (a) $x\in\Omega\setminus\overline D$ and $\sigma(]0,1])\cap\overline D=\emptyset$;
(b) $x\in\Omega\setminus\overline D$ and $\sigma(]0,1])\cap D\not=\emptyset$;
(c) $x\in\overline D$.
}\label{fig3}
\end{center}
\end{figure}

\noindent 
See Figure \ref{fig3} for an illustration of three cases.  
This theorem does not cover the case when
$x\in\Omega\setminus\overline D$ and $\sigma$ satisfies both
$\sigma(]0,1])\cap D=\emptyset$ and $\sigma(]0,1])\cap\overline
D\not=\emptyset$. However, this is quite an exceptional case.  A
similar theorem is valid in the case when $D$ is sound soft. In
Theorem 1.1 from a technical reason we needed a restriction on the
regularity of $\partial D$ ($C^2$ regularity). In Theorem 2.1 we
need only Lipschitz regularity of $\partial D$ under smallness
conditions (2.7) and (2.8) on $k$(however, being in attendance at
the competition on relaxing the regularity of $\partial D$ is not
the purpose of this paper). The piecewise linearity of the needle
is introduced just for making the geometry simple and can be
relaxed. However, from a practical point of view, it is enough.

{\it\noindent Proof.}
From Proposition 2.3 we have
$$\begin{array}{c}
\displaystyle
\int_{D}\vert v\vert^2 dy=\sum_{j}\int_{D_j}\vert v\vert^2dy\\
\\
\displaystyle
\le\sum_{j} 2\int_{D_j}\vert v-v_{A_j}\vert^2 dy
+2\int_{D_j}\vert v_{A_j}\vert^2 dy\\
\\
\displaystyle
\le\sum_{j}2C(D_j)^2(1+\frac{\vert D_j\vert^{1/2}}{\vert A_j\vert^{1/2}})^2
\int_{D_j}\vert\nabla v\vert^2 dy
+\sum_{j}2\vert D_j\vert\vert v_{A_j}\vert^2
\end{array}
$$
where $A_j\subset D_j$ and satisfy $\vert A_j\vert>0$.
\noindent
Then from Proposition 2.1 and (A.3) we have the basic inequality
$$\begin{array}{c}
\displaystyle
\int_{\partial\Omega}\{(\Lambda_0-\Lambda_{D})\overline f\}fdS
\displaystyle
\ge (1-k^2C_0(\Omega\setminus\overline D)^2)\int_{\Omega\setminus\overline D}\vert\nabla w\vert^2 dy\\
\\
\displaystyle
+\sum_{j}(1-2k^2C(D_j)^2(1+\frac{\vert D_j\vert^{1/2}}{\vert A_j\vert^{1/2}})^2)\int_{D_j}\vert\nabla v\vert^2 dy
-2k^2\vert D\vert\sum_{j}\vert v_{A_j}\vert^2.
\end{array}
\tag {2.9}
$$
Choose a sequence $\{K_n\}$ of compact sets of $\Bbb R^m$ in such a way that
$K_n\subset\Omega\setminus\sigma(]0,1])$; $\overline K_n\subset K_{n+1}$ for $n=1,\cdots$;
$\Omega\setminus\sigma(]0,1])=\cup_{n=1}^{\infty}K_n$.
Then $\vert K_n\cap D_j\vert\longrightarrow\vert D_j\setminus\sigma(]0,1])\vert
=\vert D_j\vert$ as $n\longrightarrow\infty$ uniformly with $j=1,\cdots, N$.
Thus one can take a large $n_0$ in such a way that
the set $A_j\equiv K_{n_0}\cap D_j$ satisfies
$$
\max_{j}\{2k^2(C(D_j)^2(1+\frac{\vert D_j\vert^{1/2}}{\vert A_j\vert^{1/2}})^2
-C(D_j)^2(1+1)^2)\}<\min_{j}\{1-2k^2C(D_j)^2(1+1)^2\}.
$$
We know that the sequences $\{(v_n)_{A_j}\}$ for each $j=1,\cdots, N$
are always convergent since $\overline A_j\subset\Omega\setminus\sigma(]0,1])$.
From (2.9) we have
$$
\displaystyle
I(x,\sigma,\xi)_n\ge
NC\int_{D}\vert\nabla v_n\vert^2 dy
-2k^2\vert D\vert\sum_{j}\vert (v_n)_{A_j}\vert^2
$$
where
$$
\begin{array}{c}
\displaystyle
C=\min_{j}\{1-2k^2C(D_j)^2(1+1)^2\}\\
\\
\displaystyle-
\max_{j}\{2k^2(C(D_j)^2(1+\frac{\vert D_j\vert^{1/2}}{\vert A_j\vert^{1/2}})^2-
C(D_j)^2(1+1)^2)\}>0.
\end{array}
$$
Then the blowup of $I(x,\sigma,\xi)_n$ comes from the blowup
of the sequence
$$
\displaystyle
\int_D\vert\nabla v_n\vert^2 dy.
\tag {2.10}
$$
If $x\in D$, then the blowup of the sequence given by (2.10) is
a direct consequence of Lemma 2.1.  If $x\in\partial D$, then the
exists a finite cone $V$ at vertex at $x$ such that $V\subset D$.
This is because of the Lipshitz regularity of $\partial D$. Then
Lemma 2.1 gives the blowup of the sequence.  Now consider the
case when $x\in\Omega\setminus\overline D$. If
$\sigma(]0,1])\cap\overline D=\emptyset$, then (A.3) and an
argument given in \cite{I3} provide the convergence of
$\{I(x,\sigma,\xi)_n\}$ for any needle sequence $\xi$.  If
$\sigma(]0,1])\cap D\not=\emptyset$, then there exists a point $z$
on $\sigma(]0,1[)\cap D$.  Choose an open ball centred at $z$ in
such a way that $B\subset D$.  Then from Lemma 2.2 we see the
blowup of the sequence given by (2.10).

\noindent
$\Box$

\newpage


\vspace{-1.0cm}

\noindent
As a corollary of Theorem 2.1 we obtain the characterization of $\Omega\setminus\overline D$.

\proclaim{\noindent Corollary 2.1.}  Under the same assumptions as those of Theorem 2.1
we have: a point $x\in\Omega$ belongs to $\Omega\setminus\overline D$ if and only if
there exist a needle $\sigma$ with tip at $x$ and needle sequence
$\xi$ for $(x,\sigma)$ such that
the sequence $\{I(x,\sigma,\xi)_n\}$ is bounded from above.

\endproclaim

{\it\noindent Proof.} Since we have assumed that
$\Omega\setminus\overline D$ is connected, if
$x\in\Omega\setminus\overline D$, then, one can find a piecewise
linear curve $\sigma:[0,\,1]\longrightarrow
\overline\Omega\setminus D$ with
$\sigma(0)\in\partial\Omega$, $\sigma(1)=x$ and
$\sigma(t)\in\Omega\setminus\overline D$ for all
$t\in\,]0,\,1[$. It is obvious that $\sigma$ can be chosen as an injective
curve.  This ensures the existence of a needle $\sigma$ with tip at
$x$ such that $\sigma(]0,\,1])\subset \Omega\setminus\overline D$.
Then from Theorem 2.1, one concludes the convergence of the
indicator sequence for an arbitrary needle sequence $\xi$ for
$(x,\sigma)$. Of course the existence of the needle sequence has
been ensured. If $x\in\overline D$, then again from Theorem 2.1 we
see that the indicator sequence for an arbitrary needle sequence
for $(x,\sigma)$ for an arbitrary needle $\sigma$ with tip at $x$
blows up.

\noindent
$\Box$

\section{The reflected needle-an example}

In this section we formulate a problem related to the behaviour
of the sequence of reflected solutions by an obstacle introduced below
(in the case $k=0$)
and give an explicit answer in a simple situation.
This is also an application of Lemmas 2.1 and 2.2.

{\bf\noindent Definition 3.1.} We say that the sequence
$\{g_n\}$ of $H^1(\Omega\setminus\overline D)$ functions blows up at the point $z\in\overline\Omega\setminus D$
if for any open ball $B$ centered at $z$ it holds that
$$
\displaystyle
\lim_{n\longrightarrow\infty}
\int_{B\cap(\overline\Omega\setminus D)}\vert\nabla g_n(y)\vert^2dy=\infty.
$$

\noindent
We call the set of all points $z\in\overline\Omega\setminus D$ such that
$\{g_n\}$ blows up at $z$ the blowup set of $\{g_n\}$.

\noindent
Given $x\in\Omega$, needle $\sigma$ with tip at $x$
and needle sequence $\xi=\{v_n\}$ for $(x,\sigma)$ let $u_n$ solve
$$\begin{array}{c}
\displaystyle
\triangle u=0\,\,\text{in}\,\Omega\setminus\overline D,\\
\\
\displaystyle
\frac{\partial u}{\partial\nu}=0\,\,\text{on}\,\partial D,\\
\\
\displaystyle
u=v_n\,\,\text{on}\,\partial\Omega.
\end{array}
$$
The function $u_n-v_n$ is called the reflected solution by the obstacle $D$.
It is easy to see that if $\sigma(]0,1])\cap\overline D=\emptyset$, then $\{u_n-v_n\}$
is bounded in $H^1(\Omega\setminus\overline D)$ and thus the blowup set of sequence
is empty.

We raise the following.

{\bf\noindent Problem.} What can one say about the blowup set of $\{u_n-v_n\}$
when $\sigma(]0,1])\cap\overline D\not=\emptyset$?

Here we consider the problem in a simple case in two-dimensions.
Let $R>\epsilon>0$. $\Omega$ and $D$ are given by the open discs centered at the
origin with radius $R$ and $\epsilon$, respectively. We show that,
in the case when $x\in D$, the blowup set of $\{u_n-v_n\}$ is given by
a suitable curve in $\overline\Omega\setminus D$ obtained by transforming the part
of needle $\sigma$ in $\overline D$.  We call the curve the reflected needle.

\proclaim{\noindent Proposition 3.1.} Let $\sigma$ be a needle with tip at $x\in D$
and satisfy the following: (1) $\sigma$ intersects with $\partial D$ only one time and
(2) $\displaystyle\sigma(]0,1])\cap\{y\,\vert\,\vert y\vert\le\frac{\epsilon^2}{R}\}=\emptyset$.

\begin{figure}[htbp]
\begin{center}
\epsfxsize=10cm
\epsfysize=14cm
\epsfbox{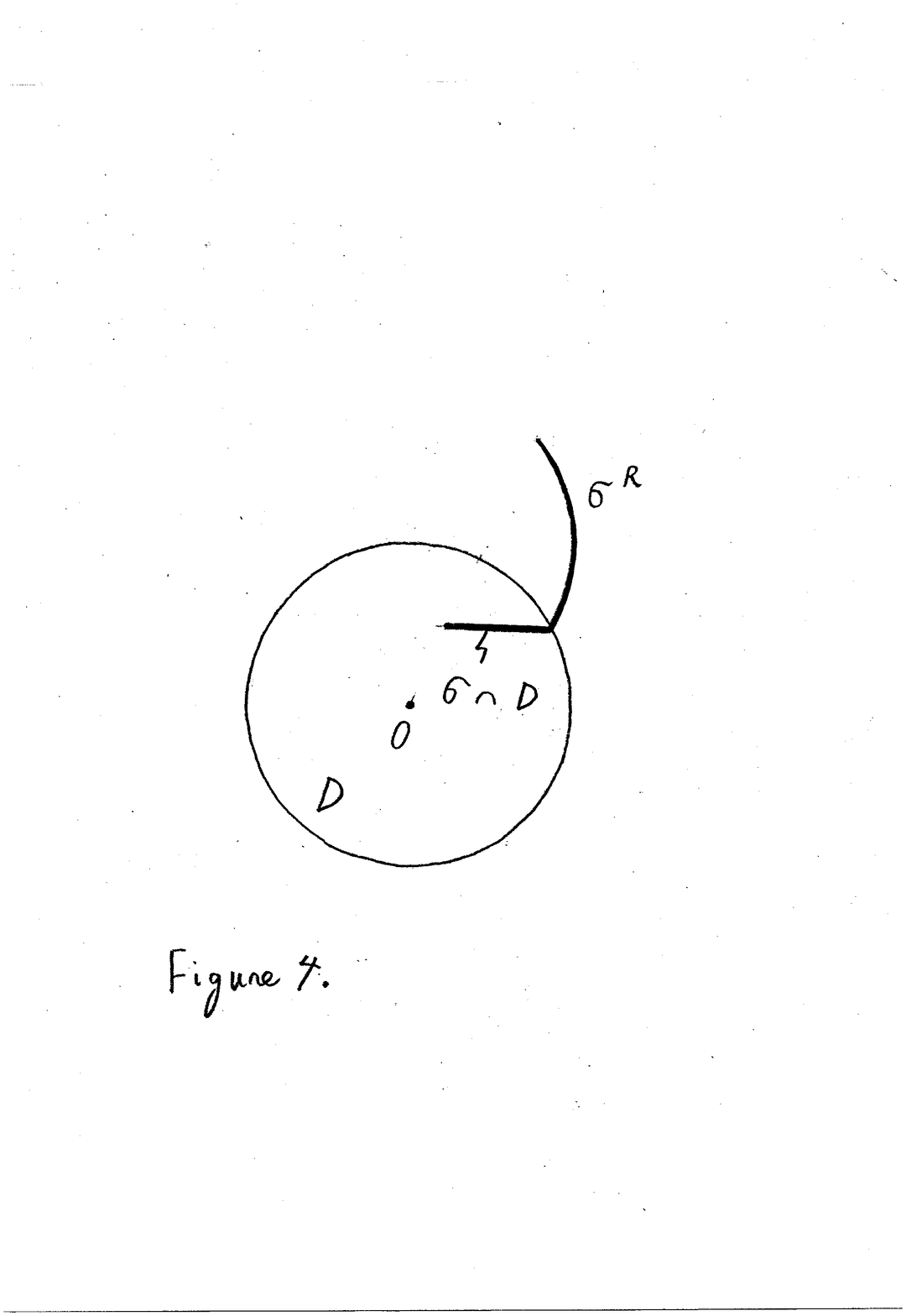}
\caption{An illustration of $\sigma^{R}$.
}\label{fig4}
\end{center}
\end{figure}

\noindent
Then the blowup set of the sequence $\{u_n-v_n\}$ coincides with the curve $\sigma^R$ given by the formula
(see Figure \ref{fig4} for an illustration of $\sigma^{R}$)
$$\displaystyle
\sigma^R=\{\frac{\epsilon^2y}{\vert y\vert^2}\,\vert\,y\in\sigma(]0,1])\cap\overline D\}.
$$
\endproclaim

{\it\noindent Proof.} Choose $\varphi\in C^{\infty}_0(\Bbb R^2)$
in such a way that $\varphi=1$ in a neighbourhood of
$\sigma(]0,1])\cap\overline D$ and $\varphi=0$ in a neighbourhood
of the circle centered at the origin with radius $\displaystyle\epsilon^2/R$.
\noindent
Given needle sequence $\xi=\{v_n\}$ for $(x,\sigma)$ define
$$
\displaystyle
\tilde{v_n}(z)=\varphi(y)v_n(y),\,\,z\in\Omega\setminus\overline D
$$
where
$$
y=\frac{\epsilon^2 z}{\vert z\vert^2}.
$$
Note that this is nothing but the Kelvin transform of the function $\varphi v_n$
with respect to the circle centered at the origin with radius $\epsilon$.

Set
$$
R_n(z)=u_n(z)-v_n(z)-\tilde{v_n}(z),\,\,z\in\Omega\setminus\overline D.
$$
This function vanishes on $\partial\Omega$. A direct computation
by using the polar coordinates around the origin gives the formula
$$\displaystyle
\triangle_z\{g(y)\}
=\frac{\epsilon^{4}}{\vert z\vert^{4}}
(\triangle g)(y)
$$
where $g$ is an arbitrary function in $\Bbb R^2\setminus\{0\}$.
Applying this formula to $g=\varphi v_n$, we have
$$\displaystyle
(\triangle R_n)(z)=-\frac{\epsilon^{4}}{\vert z\vert^{4}}\{
(\triangle\varphi)(y)v_n(y)+2\nabla\varphi(y)\cdot\nabla v_n(y)\}.
$$
This right-hand side is convergent as $n\longrightarrow\infty$ since both $\nabla\varphi(y)$
and $\triangle\varphi(y)$ vanish in a neighbourhood of the curve $\sigma(]0,1])\cap\overline D$;
both $v_n(y)$ and $\nabla v_n(y)$ are convergent in $L^2(K)$ as $n\longrightarrow\infty$
where
$$\displaystyle
K=\{y\,\vert\,\frac{\epsilon^2}{R}\le \vert y\vert\le\epsilon\,\,\text{and}\,\,
\text{dis}\,(y,\sigma(]0,1])\cap\overline D)\ge\eta\}\subset\Omega\setminus\sigma(]0,1])
$$
and $0<\eta$.

A direct computation also gives
$$\displaystyle
\frac{\partial R_n}{\partial\nu}(z)=-(1-\varphi(z))\frac{\partial v_n}{\partial\nu}(z)
+\frac{\partial\varphi}{\partial\nu}(z)v_n(z),\,\,\vert
z\vert=\epsilon.
$$
This right hand sid is convergent in $H^{-1/2}(\partial D)$ since
both $\partial\varphi/\partial\nu$ and $1-\varphi$ vanish for $z$
close to the single point in the set $\sigma(]0,1])\cap\partial
D$.  Then the well posedness of the mixed boundary value problem
yields that the sequence $\{R_n\}$ is bounded in
$H^1(\Omega\setminus\overline D)$. Then from Lemmas 2.1 and 2.2
one obtains the desired conclusion.

\noindent
$\Box$

\newpage


\vspace{-1.0cm}

\noindent
We think that Proposition 3.1 is a special case of a more general theorem
that shall give the description of the blowup set of $\{u_n-v_n\}$
by using a curve obtained by {\bf a rule}.
In a forthcoming paper we will consider the problem of seeking such a rule for general
$D$, $\Omega$ and $k>0$.

\section{Remark}

\noindent
It is possible to obtain the corresponding results in other applications
of the probe method (see \cite{I1, I4}).  For example, consider the Dirichlet-to-Neumann
map $\Lambda_{\gamma}$ for the equation $\nabla\cdot\gamma\nabla u=0$ in $\Omega$. Here $\gamma=\gamma(y)$
denotes the electrical conductivity.  Assume that $\gamma$ takes the form:
$\gamma(y)=1,\,y\in\,\Omega\setminus D$ and $1+h(y),\,\,y\in\,D$ where $h(y)$ is given by
a function in $L^{\infty}(D)$ satisfying $\text{ess}\,\inf_{y\in D}(1+h(y))>0$ and
the global jump condition: $h(y)\ge C$ a.e. in $D$ or $-h(y)\ge C$ a.e. in $D$ for a positive constant $C$.
We obtain
\proclaim{\noindent Theorem 4.1.}
A point $x\in\Omega$ belongs to $\Omega\setminus\overline D$ if and only if
there exist a needle $\sigma$ with tip at $x$ and
needle sequence $\xi=\{v_n\}$ for $k=0$ such that
the sequence $\{I(x,\sigma,\xi)_n\}$ given by the formula
$$\begin{array}{c}
\displaystyle
I(x,\sigma,\xi)_n=\int_{\partial\Omega}\{(\Lambda_{\gamma}-\Lambda_1)f_n\}f_ndS,\\
\\
\displaystyle
f_n(y)=v_n(y),\,\,y\in\partial\Omega,
\end{array}
$$
is bounded.  Moreover given $x\in\Omega$ and needle $\sigma$ with tip at $x$ we have that

if $x\in\Omega\setminus\overline D$ and $\sigma(]0,1])\cap\overline D=\emptyset$,
then for any needle sequence $\xi=\{v_n\}$ for $(x,\sigma)$
the sequence $\{I(x,\sigma,\xi)_n\}$ is convergent

if $x\in\Omega\setminus\overline D$ and $\sigma(]0,1])\cap D\not=\emptyset$,
then for any needle sequence $\xi=\{v_n\}$ for $(x,\sigma)$ we have
$\lim_{n\longrightarrow\infty}\vert I(x,\sigma,\xi)_n\vert=\infty$

if $x\in\overline D$,
then for any needle sequence $\xi=\{v_n\}$ for $(x,\sigma)$ we have
$\lim_{n\longrightarrow\infty}\vert I(x,\sigma,\xi)_n\vert=\infty$.

\endproclaim

\noindent

This theorem suggests that the new formulation of the probe method
can probably be considered as a final generalization of the enclosure method
introduced in \cite{IE}.
The needle sequences play the role similar to the special harmonic functions
coming from Mittag-Leffler's function in a generalized enclosure method
given in \cite{I5, IS}.  The proof is a direct consequence of the system of the integral
inequalities (\cite{I}) and Lemmas 2.1 and 2.2 for $k=0$.

We point out that the behaviour of $\{I(x,\sigma,\xi)_n\}$ for general $x\in \overline D$ is not clear
without a global assumption on $h$ in $D$.
However, one can easily deduce that if $h$ or $-h$ has a
positive lower bound in a neighbourhood of $\sigma(]0,\,1])\cap
\overline D$, then $\lim_{n\longrightarrow\infty}\vert I(x,\sigma,\xi)_n\vert=\infty$.

In my opinion, it is impossible to
know the behaviour of $I(x,\sigma,\xi)_n$ as $n\longrightarrow\infty$
for $x\in\overline D$ from the property of the needle sequence in the case when both $h$ and $-h$ do not have a
positive lower bound in any neighbourhood of  $\sigma(]0,\,1])\cap
\overline D$. For this purpose we have to study the behaviour of the
sequence of reflected/refracted solutions
by the obstacles, inclusions and cracks.  We also think that {\bf the study may enable us to
drop the restriction on $k$ given by (2.7) and (2.8).}

$$\quad$$

\centerline{{\bf Acknowledgement}}

This research was partially supported by Grant-in-Aid for Scientific
Research (C)(2) (No.  15540154) of Japan  Society for the Promotion of Science.
The author thanks Klaus Erhard and Roland Potthast for providing me with a preprint of \cite{EP}
and the anonymous referees for their valuable comments and suggestions for improvement
of the original manuscript.

$$\quad$$

\section{Appendix}

{\bf\noindent A.1.  Remark.}
In the proof of Theorem 4 of \cite{I3} some
important explanations described below are missing.

(1) $f$ in (A.1) of the paper should belong to $\{H^1(U)\}^*$ and satisfy
$f(v\vert_U)=0$ for all $v\vert_U\in\,Y$;

(2) the right hand side of (A.1) of the paper defines a bounded linear functional on $H^1_0(\Omega)$.

{\bf\noindent A.2.  Estimates.}

\proclaim{\noindent Proposition A.1.}
Let $W$ be a bounded domain with $C^2$ boundary of $\Bbb R^m$.
Let $v\in H^1(W)$ satisfy $\triangle v+k^2v=0$ in $W$.
Then, for any compact set $K$ of $\Bbb R^m$ with $K\subset W$ there
exists a positive constant $C'=C'(K,W)$ independent of $v$ such that
$$\displaystyle
\int_K\vert\nabla v\vert^2dy
\le C'\int_W\vert v\vert^2 dy.
\tag {A.1}
$$
Moreover, assume that $0$ is not a Dirichlet eigenvalue of $\triangle+k^2$ in $W$.
Then there exists a positive constant $C=C(W)$ independent of $v$
such that
$$\displaystyle
\int_W\vert v\vert^2 dy\le C\int_{\partial W}\vert v\vert^2 dS.
\tag {A.2}
$$

\endproclaim

{\it\noindent Proof.}
First we prove (A.1).  Let $\varphi\in C^{\infty}_0(W,\Bbb R)$.  Multiply the equation
$\triangle v+k^2v=0$ in $W$ by $\varphi^2\overline v$ and integrate the resultant
equation over $W$.  Integration by parts gives
$$
\displaystyle
\int_W\vert\nabla v\vert^2\varphi^2dy
\le C\int_W\vert v\vert^2\vert\nabla\varphi\vert^2 dy
+\int_W k^2\vert v\vert^2\varphi^2dy.
$$
Choose $\varphi\in C^{\infty}_0(W)$ in such a way that
$\varphi=1$ on $K$ and $0\le\varphi\le 1$.  Then one gets (A.1).

Let $z\in H^2(W)$ solve
$$\begin{array}{c}
\displaystyle
\triangle z+k^2z=\overline v,\,\,\text{in}\,W,\\
\\
\displaystyle
z=0\,\,\text{on}\,\partial W.
\end{array}
$$
Then we have
$$\begin{array}{c}
\displaystyle
\int_W\vert v\vert^2 dy=\int_W(\triangle z+k^2z)vdy\\
\\
\displaystyle
=\int_{\partial W}\frac{\partial z}{\partial\nu}vdS
-\int_W\nabla z\cdot\nabla vdy
+\int_W k^2zv\\
\\
\displaystyle
=\int_{\partial W}\frac{\partial z}{\partial\nu}vdS
\end{array}
$$
and the trace theorem yields
$$\begin{array}{c}
\displaystyle
\Vert v\Vert^2_{L^2(W)}\le\int_{\partial W}\vert\frac{\partial z}{\partial\nu}\vert\vert v\vert dS\\
\\
\displaystyle
\le \Vert\nabla z\Vert_{L^2(\partial W)}\Vert v\vert_{\partial W}\Vert_{L^2(\partial W)}\\
\\
\displaystyle
\le C_1\Vert z\Vert_{H^2(W)}\Vert v\vert_{\partial W}\Vert_{L^2(\partial W)}\\
\\
\displaystyle
\le C_2\Vert v\Vert_{L^2(W)}\Vert v\vert_{\partial W}\Vert_{L^2(\partial W)}.
\end{array}
$$
Thus we obtain (A.2).

\noindent
$\Box$

The reader can see this type of argument for the proof of this proposition, e.g.,
in \cite{KV}.

{\bf\noindent A.3.  An integral identity.}

\noindent
The identity below has been established in \cite{I3}.

\proclaim{\noindent Proposition A.2.}
For all $f\in H^{1/2}(\partial\Omega)$
$$\begin{array}{c}
\displaystyle
\int_{\partial\Omega}\{(\Lambda_0-\Lambda_{D})\overline f\}fdS\\
\\
\displaystyle
=\int_{\Omega\setminus\overline D}\vert\nabla(u-v)\vert^2dy
-k^2\int_{\Omega\setminus\overline D}\vert u-v\vert^2 dy\\
\\
\displaystyle
+\int_{D}\vert\nabla v\vert^2 dy-k^2\int_{D}\vert v\vert^2 dy
\end{array}
\tag {A.3}
$$
where $u$ solves
$$\begin{array}{c}
\displaystyle
(\triangle +k^2)u=0\,\,\text{in}\,\,\Omega\setminus\overline D,\\
\\
\displaystyle
u=f\,\,\text{on}\,\,\partial\Omega,\\
\\
\displaystyle
\frac{\partial u}{\partial\nu}=0\,\,\text{on}\,\,\partial D;
\end{array}
$$
$v$ solves
$$\begin{array}{c}
\displaystyle
(\triangle +k^2)v=0\,\,\text{in}\,\,\Omega,\\
\\
\displaystyle
v=f\,\,\text{on}\,\,\partial\Omega.
\end{array}
$$
\endproclaim

\vskip1cm
\noindent
e-mail address

ikehata@math.sci.gunma-u.ac.jp
\end{document}